\documentclass[a4paper,12pt]{article}
\usepackage[T2A]{fontenc}
\usepackage[english,russian]{babel}
\usepackage{amsfonts,amssymb,amscd,amsmath}
\usepackage[dvips]{graphicx}
\usepackage[right=25mm,left=25mm,top=30mm,
bottom=40mm]{geometry}

\begin{document}
\begin{center}
{\Large \bf Теоремы типа Бохера
для динамических уравнений на
временных шкалах}
\end{center}
%\begin{center}
%{\Large \bf Theorems of B\^ocher's 
%type for dynamic equations on 
%time scales}
%\end{center}

\vskip 12pt\noindent

\centerline{\bf \copyright 
Владимир Шепселевич Бурд}

\vskip 12pt\noindent

\centerline{Ярославский государственный
университет им. П.Г. Демидова}

\vskip 14pt\noindent

{\bf Аннотация.} Найдены условия,
при которых все решения систем
динамических уравнений на 
временных шкалах стремятся к 
конечным пределам при 
$t\to\infty$.

\vskip 14pt\noindent

\begin{center}
{\Large \bf Theorems of B\^ocher's 
type for dynamic equations on 
time scales}
\end{center}

\centerline{\bf V.Sh.~Burd}

\vskip 12pt\noindent

\centerline{Demidov Yaroslavl State University}

\vskip 12pt\noindent

\bigskip

{\bf Abstract}\newline
The conditions are found that
all solutions of a systems 
dynamic equations on time scales
tends to finite limits as 
$t\to\infty$. 
%\end{quotation}

\vskip 12pt\noindent

A. Wintner [1,2] has named by
theorems of type B\^oher
a assertions which is guaranteed
that the each nontrivial solution
of the linear system of ordinary
differential equations
$$
\frac{dx}{dt}=A(t)x, \quad
t_0\le t<\infty,
$$
where $A(t)$ is $n\times n$ 
matrix, tends to nontrivial 
finite limit as $t\to\infty$.
It has place, if the elements 
of matrix $A(t)$ belong to
$L_1[t_0,\infty]$. This result
is connected with name of B\^ocher.
[3] (see also [4]). А. Wintner
has receive series of less 
restrictive conditions.
P, Hartman [5] has establish
analogue of B\^оher's
result for nonlinear system.
\par
In this paper analogical problem
is studied for dynamical equations
on time scales. 
\par
We consider the linear equation
$$
x^\Delta=A(t)x,\quad t\in T,
\eqno(1)
$$
where $T\subset{\cal R}$ is 
time scale and
$\sup T=\infty$, $A(t)$ be an
$n\times n$-valued function on 
$T$. It is assumed that $A(t)$ is
rd-continuous and regressive
(see [6]).
\par
The system(1) will be said to be
of class (S) if (i) every 
solution $x(t)$ of (1) has a 
limit $x_\infty$ as $t\to\infty$,
and (ii) for every constant 
vector $x_\infty$ there is a 
solution $x(t)$ of (1)
such that $x(t)\to x_\infty$ as 
$t\to\infty$.
\par 
Evidently that (1) is of class 
(S) if and only if for every 
fundamental matrix $X(t)$ of (1),
$X_\infty=\lim X(t)$ exists as 
$t\to\infty$ and is nonsingular.
\par
Let's note still that (1) is of 
class (S), if it has a 
fundamental matrix of form
$$
X_0(t)=I+o(1)
$$
as $t\to\infty$.
\par
Let $a\in T,\,a>0$. Let integral
$\int\limits_a^\infty A(s)ds$ is
convergent absolutely. Following 
theorem well known.

\bigskip

{\bf Theorem 1}. {\it If integral
$\int\limits_a^\infty A(s)\Delta 
s$ is convergent absolutely. 
Then system (1) is of class 
(S)}.

\bigskip

{\bf Proof}. The solution of 
system (1) are represented in 
the form
$$                                    
x(t)=x(a)+\int\limits_a^tA(s)
x(s)\Delta s. \eqno(2) 
$$                                    
We have the inequality
$$
|x(t)|\le|x(a)|+\int\limits_a^t
|A(s)||x(s)|\Delta s. \eqno(3) 
$$
The Gronwall's inequality for
times scales gives  
$$                                                    
|x(t)|\le|x(0)|e_{|A|}(t,a),
\quad t\in T,
$$           
where $e_{|A|}(t,a)$ is 
exponential function of 
equation
$$
x^\Delta=|A(t)|x.
$$                
Consequently all solutions of 
system (1) is bounded.                                            
Let $t_2>t_1$. Let's evaluate 
norm $|x(t_2)-x(t_1)|$:                                
$$                                                                           
|x(t_2)-x(t_1)|\le\int\limits_
{t_1}^{t_2}|A(s)||x(s)|\Delta s. 
\eqno(4)            
$$                                                                           
From inequality (4) follows that
$|x(t_2)-x(t_1)|$ can be made as 
much as small if to take $t_1$ 
big enough,
Therefore $|x(t)|$ converges to 
the finite limits as $t\to\infty$.
From Gronwall's inequality 
follows
$$
|x(t)|\ge |x(s)|
(e_{|A|}(t,s))^{-1},\quad t,s
\in T. \eqno(5)
$$
The inequality (5) implies that 
$x_\infty\ne 0$ unless 
$x(n)\equiv 0$. This proves
theorem 1.

\bigskip

Let now integral $\int\limits_
a^\infty A(s)\Delta s$ is 
convergent
(possibly just conditionally).
In system (1) make the change of 
variables
$$
x(t)=y(t)+Y(t)y(t),
$$
where $N\times N$ matrix $Y(t)$
we shall choose later. We obtain
$$
[I+Y(\sigma(t))]y^\Delta+
Y^\Delta(t)y(t)=A(t)y(t)+
A(t)Y(t)y(t).   \eqno(6)
$$
Let
$$
Y^\Delta(t)=A(t).
$$
Then
$$
Y(t)=-\int\limits_t^\infty A(s)
\Delta s.
$$
Evidently $Y(t)\to 0$ as 
$t\to\infty$. Therefore matrix 
$I+Y(\sigma(t))$ has a bounded
converse for sufficiently large 
$t$. Therefore
we obtain for sufficiently large 
$t$
$$
y^\Delta(t)=(I+Y(t))^{-1}A(t)
Y(t)y(t). \eqno(7)
$$

\bigskip

{\bf Theorem 2}. {\it Let integral
$\int\limits_a^\infty A(s)
\Delta s$ is 
convergent and integral 
$$
\int\limits_a^\infty A(t)
(\int\limits_{s=t}^\infty
 A(s)\Delta s)\Delta t    \eqno(8)
$$
is convergent absolutely. Then 
system (1) is of class (S).}

\bigskip

{\bf Proof}. From conditions 
theorem 2
follows that right part of system
(7) satisfy conditions of 
theorem 1.

\bigskip

{\bf Remark.} If $X(t)$ is 
fundamental matrix of system (1),
then $X^{-1}(t)$ is fundamental 
matrix of the system
$$
x^\Delta=x'(\ominus A)(t),
$$
where (see [6])
$$
(\ominus A)(t)=-A(t)
[I+\mu(t)A(t)]^{-1}.
$$
Therefore the theorem 2 remains
fair, if the requirement of 
absolute convergence integral 
(8)
to replace with the requirement
of absolute convergence of 
integral 
$$
\int\limits_a^\infty 
(\int\limits_{s=t}^\infty
(\ominus A)(s)\Delta s)A(t)
\Delta t.
$$
\par
Further we receive the theorem.

\bigskip

{\bf Theorem 3}. {\it Let integrals
$\int\limits_a^\infty A(s)
\Delta s$  
and 
$\int\limits_a^\infty A(t)
(\int\limits_{s=t}^\infty
 A(s)\Delta s)\Delta t$ are 
convergent and integral
$$
\int\limits_a^\infty A(t)
(\int\limits_{s=t}
^\infty A(s)(\int\limits_{\tau=s}
^\infty A(\tau)\Delta\tau)\Delta
 s)\Delta t
$$
is convergent absolutely. 
Then system (1) is of class (S).}

\bigskip

{\bf Proof}. In system (1) make 
now the change of variables
$$
x(t)=y(t)+Y_1(t)y(t)+Y_2(t)y(t).
$$
where matrices $Y_1(t),\,Y_2(t)$ are
we shall choose later. We obtain
$$
\begin{array}{l}
(I+Y_1(\sigma(t))+Y_2(\sigma(t))
y(t)^\Delta +
(Y_1^\Delta(t))y(t)+\\+
(Y_2^\Delta(t))y(t)=
(A(t)+A(t)Y_1(t)+A(t)Y_2(t))y(t)
\end{array}. \eqno(9)
$$
Let
$$
Y_1^\Delta(t)=A(t),\quad 
Y_2^\Delta(t)=A(t)Y_1(t).
$$
Then
$$
Y_1(t)=-
\int\limits_t^\infty A(s)
\Delta s.
\quad Y_2(t)=
\int\limits_t^\infty A(s)
(\int\limits_{\tau=s}^\infty
A(\tau)\Delta\tau)\Delta s 
\eqno(10) 
$$
From the formulas (10) follows 
that $Y_1(\sigma(t)),\,
Y_2(\sigma(t))$ are converges to
0 as $t\to\infty$.
Therefore matrices
$I+Y_1(\sigma(t)),\,I+
Y_2(\sigma(t))$ have a bounded
converse for sufficiently large
$t$.                                                                                                             
Hence system (9) can be written in
following form
$$
y^\Delta(t)=(I+Y_1(\sigma(t))+
Y_2(\sigma(t)))^
{-1}A(t)Y_2(t)y(t). \eqno(11)
$$
From conditions theorem 3
follows that right part of system
(6) satisfy conditions of 
theorem 1.

\bigskip

Generally, if integrals
$$
\int\limits_a^\infty A(t)\Delta 
t,\quad 
\int\limits_a^\infty A(t)(\int\limits
_{s=t}^\infty A(s)\Delta s)
\Delta t,\dots,
$$
$$
\int\limits_a^\infty A(t_1)
(\int\limits_{t_2=t_1}^\infty 
A(t_2)\Delta t_2
\cdots
\int\limits_{t_k=t_{k-1}}^\infty
A(t_k)\Delta t_k)\Delta t_1
$$
are convergent and integral
$$
\int\limits_a^\infty A(t_1)
\int\limits_{t_2=t_1}^\infty A(t_2)
\Delta t_2
\cdots\int\limits_{t_{k+1}=
t_k}^\infty A(t_{k+1})\Delta
t_{k+1})\Delta t_1 
$$
is convergent absolutely, then 
system (1) is of class (S).
\par
Consider now the nonlinear system of difference equations
$$
x^\Delta(t)=f(t,x(t)), \quad 
t\in T,\, x\in R^N.  
\eqno(12)
$$

\bigskip

{\bf Theorem 4}. {\it Let $f(t,
x)$ be defined for $t\in T\,|x|
<\delta\, (\le\infty)$ and 
satisfy inequality
$$
|f(t,x)|\le K(t)|x|,
$$
where
$$
\int\limits_a^\infty K(t)\Delta t
<\infty.
$$
If $|x_0|$ is sufficiently 
small, let us
say
$$
|x_0|e_K(t,t_0)<
\delta, \eqno(13)
$$
Then for a solution $x(t)$ of (12) 
satisfying $x(a)=x_0$ exists
$$
x_\infty=\lim_{t\to\infty}x(t)
$$
and $x_\infty\ne 0$ unless $x(t)
\equiv 0$.}

\bigskip

{\bf Proof}. The solutions of 
system (12) are represented in 
the form
$$
x(t)=x(a)+\int\limits_a^t
f(s,x(s))\Delta s.  \eqno(14)
$$
If $x(a)=x_0$ satisfies (14),
then from (13) and conditions of 
theorem 4 follows that
$$
|x(t)|\le[|x(a)|+\int\limits_a^t
K(t)|x(t)|.
$$
The Gronwall's inequality gives
$$
|x(t)|\le|x(a)e_K(t,t_0).
$$
Consequently all solutions of 
system (11) is bounded.
Let $t_2>t_1$. Let's evaluate 
norm $|x(t_2)-x(t_1)|$:                                
$$                                                                           
|x(t_2)-x(t_1)|\le\int\limits_
{t_1}^{t_2}K(s)|x(s)|\Delta s. 
\eqno(15)            
$$                                                                           
From inequality (4) follows that
$|x(t_2)-x(t_1)|$ can be made as 
much as small if to take $t_1$ 
big enough,
Therefore $|x(t)|$ converges to 
the finite limits as $t\to\infty$.
From Gronwall's inequality 
follows
$$
|x(t)|\ge |x(s)|
(e_K(t,s))^{-1},\quad t,s
\in T. \eqno(16)
$$
The inequality (16) implies that 
$x_\infty\ne 0$ unless 
$x(n)\equiv 0$. This proves
theorem 4.

\bigskip

Let $T={\cal R}$. Then the 
theorem 2 together with the
remark 1 return to the results 
Wintner [1].

\bigskip

Let $T={\cal Z}$. We consider
a discrete adiabatic oscillator
$$
x(n+2)-(2\cos\alpha)x(n+1)+
(1+g(n))x(n)=0, \eqno(17)
$$
where $0<\alpha<\pi,\,g(n)\to 0$
as $n\to\infty$. For 
$g(n)\equiv 0$ the equation (16)
has the form
$$
x(n)=C_1\cos n\alpha+C_2\sin n
\alpha.
$$
We convert (16) into a system of
equations by introducing new 
variables $C_1(n),\,C_2(n)$
$$
x(n)=C_1(n)\cos n\alpha+C_2(n)
\sin n\alpha,
$$
$$
x(n+1)=C_1(n)\cos(n+1)\alpha+
C_2(n)\sin(n+1)\alpha.
$$
We obtain the system
$$
\Delta u(n)=B(n)g(n)u(n),
$$ 
where $\Delta u_1(n)=C_1(n+1)-
C_1(n),\,\Delta u_2(n)=C_2(n+1)-
C_2(n)$. The matrix $B(n)$ has 
the form
$$
B(n)=\frac{1}{\sin \alpha}(A_0+
A_1(n)),
$$
where 
$$
A_0=\left(\begin{array}{cc}
\sin\alpha&
\cos\alpha\\-\cos\alpha&\sin\alpha
\end{array}\right),\quad
A_1(n)=\left(\begin{array}{cc}
\sin(2n+1)\alpha&
-\cos(2n+1)\alpha\\-\cos(2n+1)
\alpha&-\sin(2n+1)\alpha\end{array}
\right).
$$

\bigskip

{\bf References}

%\begin{enumerate}
%\item
1. Wintner A. On a theorem of 
B\^ocher in the theory of ordinary
linear  differential equations,
Amer. J. Math., v. 76,
1954, pp. 183 - 190.
%\item

2. Wintner A. Addenda to the paper 
on B\^ocher's theorem, Amer. J. 
Math., v. 78, 1956, pp. 895 - 897.
%\item

3. B\^ocher M. On regular singular
points of  linear differential
equations of the second order 
whose coefficients are not 
necessarily  analytic, Trans. 
Amer. Math. Soc., v.1, 1900, pp.
40 - 53.
%\item

4. Dunkel O. Regular singular 
points of a system  of 
homogeneous linear differential 
equations of the first order,
Proc. Amer. Acad.  Arts sci., v.
38, 1902-3, pp. 341 - 370. 
%\item

5. Hartman P. Ordinary differential 
equations,New York, Wiley. 1964,
612 p. 
%\item

6. Bohner M., Peterson A. Dynamic
equations on time scales, Boston,
Birkhauser, 2000, 358 p.
%\end{enumerate}

\end{document}